\documentstyle[11pt]{article}

\def\el {{\`e}}
\def\ep {{\'e}}

\begin{document}

\noindent{\small \bf Syst\el mes dynamiques/Dynamical systems
\\
\noindent \'Equations diff\ep rentielles/Ordinary Differential
Equations}

 \vspace{0.2cm}

\noindent{\LARGE \bf  La non--int\ep grabilit\'e m\ep romorphe du
probl\el me plan des trois corps} \vspace{0.3cm}

\noindent {\Large \bf  Tsygvintsev Alexei}

\vspace{0.2cm}

{\small \bf \noindent Laboratoire de Math\ep matiques Emile
Picard, C.N.R.S.-- UMR 5580, Universit\'e Paul Sabatier,  118,
route de Narbonne, 31062 Toulouse C\ep dex 4, France \\ \noindent
Courriel : tsygvin@picard.ups-tlse.fr \\

\noindent  et Section de Math\ep matiques, Universit\'e de Gen\el
ve, 2-4, rue du Li\el vre, Case postale 240, CH-1211 Gen\el ve 24,
Suisse  \\ \noindent Courriel : Alexei.Tsygvintsev@math.unige.ch }

\vspace {0.7cm}

 \noindent{\bf R\ep sum\'e.} {\small Nous consid\ep
rons le probl\el me des trois corps dans le plan et montrons que
ce probl\el me ne poss\el de pas  de syst\el me complet d'int\ep
grales premi\el res  m\ep romorphes \`a un voisinage de la
solution particuli\el re de Lagrange.}

\vspace{0.2cm} \noindent { \bf The meromorphic non--integrability
of the planar three--body problem } \vspace{0.2cm}

\noindent{\bf Abstract.}{\small We study the planar three-body
problem and prove the absence of a complete set of complex
meromorphic first integrals in a neighborhood of the Lagrangian
solution.}

\vspace{0.3cm}

 \noindent { \bf 1. Introduction et  r\ep sultats }

 Consid\ep rons  trois corps $P_1$, $P_2$, $P_3$  sur le  plan
avec des masses $m_1$, $m_2$, $m_3$ qui s'attirent conform\ep ment
\`a la loi de Newton. Soient
 $(x_1,x_2)$ les coordonn\ep es de $P_1$, $(x_3,x_4)$ les
coordonn\ep es de $P_2$, $(x_5,x_6)$ les coordonn\ep es de $P_3$.
Notons $y_r=\displaystyle m_k\frac{dx_r}{dt}$ o\`u $k$ est la
partie enti\el re de  $(r+1)/2$.

Les \ep quations de Hamilton du mouvement des corps prennent  la
forme $$ \displaystyle\frac{dx_r}{dt}=\displaystyle\frac{\partial
H_1}{\partial y_r}, \quad
\displaystyle\frac{dy_r}{dt}=-\displaystyle\frac{\partial
H_1}{\partial x_r}, \quad (r=1,2,\dots,6), \eqno (1.1)  $$ avec
l'Hamiltonien $$
\begin{array}{ll}
H_1=\displaystyle
  \frac{1}{2m_1}(y_1^2+y_2^2)+\displaystyle \frac{1}{2m_2}(y_3^2+y_4^2)+\displaystyle \frac{1}{2m_3}(y_5^2+y_6^2)- m_3m_2\{(x_3-x_5)^2+
(x_4-x_6)^2\}^{-1/2}\\-m_3m_1\{(x_5-x_1)^2+(x_6-x_2)^2\}^{-1/2}-
m_1m_2\{(x_1-x_3)^2+ (x_2-x_4)^2\}^{-1/2}.
\end{array}
$$

Ce probl\el me peut \^etre illustr\'e par l'attraction mutuelle de
la terre, de la lune et du soleil.

Les int\ep grales premi\el res connues du syst\el me (1.1)  sont

\noindent $ F_1=H_1$ -- l'\'energie,\\ $ F_2=y_1+y_3+y_5,$ $
F_3=y_2+y_4+y_6$ -- les int\ep grales du mouvement du centre de
gravit\' e,\\ $ F_4=y_1x_2+y_3x_4+y_5x_6-x_1y_2-x_3y_4-x_5y_6=k$
-- l'int\ep grale des aires.

En supposant que le centre de gravit\'e est fixe, c'est-\`a-dire
que $F_2=F_3=0$, et en utilisant l'int\ep grale des aires $F_4$,
le syst\el me (1.1) peut \^etre ramen\'e \`a $3$ degr\ep s  de
libert\'e [6] $$
\displaystyle\frac{dq_r}{dt}=\displaystyle\frac{\partial
H}{\partial p_r},\quad \displaystyle\frac{dp_r}{dt}=
-\displaystyle\frac{
\partial H}{\partial q_r}, \quad (r=1,2,3), \eqno (1.2)
$$ o\`u
 $$
\begin{array}{ll}
H=\displaystyle \frac{M_1}{2}\left\{p_1^2+\displaystyle
\frac{1}{q^2_1}P^2\right\}+\displaystyle \frac{M_2}{2}(p_2^2+
p_3^2)+\displaystyle \frac{1}{m_3}\left\{p_1p_2-\displaystyle
\frac{p_3}{q_1}P\right\} -\displaystyle \frac{ m_1m_3}{r_1}-
\displaystyle \frac{m_3m_2}{r_2}- \displaystyle \frac{
m_1m_2}{r_3}, \\ P=p_3q_2-p_2q_3-k,
\end{array}
$$  et $$ r_1=q_1, \quad r_2=\sqrt{q^2_2+q^2_3}, \quad r_3=\sqrt{
(q_1-q_2)^2+q^2_3},$$ sont les distances mutuelles des corps.

Nous appellerons le syst\el me (1.2) le {\it probl\el me plan des
trois corps}.

D'apr\`es  Bruns [1] le probl\el me (1.1) et donc le probl\el me
(1.2) n'admet pas d'int\ep grale premi\el re alg\ep brique, autre
que  les int\ep grales premi\el res d\ep ja connues.

Poincar\'e [4] a d\ep montr\'e en 1892  que si la troisi\el me
masse $m_3=\mu$ est infiniment petite et est attire\' e par les
deux premi\el res masses $m_1$  et $m_2$, alors  le probl\el me
(1.2) n'admet  pas d'int\ep grale premi\el re holomorphe en $\mu$
autre que celle de l'\ep nergie.

Il est bien connu  (Lagrange [2]) que ce probl\el me  poss\el de
une solution particuli\el re dans laquelle  les trois corps
forment un triangle \ep quilat\ep ral et d\'ecrivent chacun une
conique. Notons $\Gamma$ celle qui correspond au cas parabolique.

Nous consid\ep rons dans cette Note la question de l'existence de
deux  int\'egrales premi\el res supplementaires  m\ep romorphes et
fonctionnellement ind\ep pendantes autres que  celle de l'\ep
nergie au voisinage de la solution $\Gamma$.

\noindent Notre r\ep sultat principal est le suivant:\\

 \noindent{\it TH\'EOR\`EME 1.1--}  Le probl\el me
plan des trois corps  n'admet pas  deux int\ep grales premi\el res
supplementaires m\ep romorphes  par rapport aux variables $q_i$,
$p_j$, $r_k$ et fonctionnellement ind\ep pendantes  au voisinage
de la solution particuli\el re de Lagrange $\Gamma$.

\vspace{0.5cm}

\noindent{\it COROLLAIRE 1.2--} Le probl\el me plan des trois
corps  n'est pas compl\el tement m\ep romorphiquement   int\ep
grable au voisinage de la solution $\Gamma$.

\vspace{0.5cm}

La preuve est bas\ep e sur la th\ep orie  de Ziglin [7] qui donne
une condition n\ep cessaire de non--int\ep grabilit\' e.

Dans le cas de masses \ep gales, le m\^eme r\ep sultat a \ep t\'e
obtenu par D. Boucher  \`a l'aide du th\ep or\el me de non--int\ep
grabilit\' e de Moralis-Ramis [3].

\vspace{0.5cm}

\noindent { \bf 2. Les \ep quations normales variationnelles  le
long de la solution $\Gamma$ }

\vspace{0.5cm}

 Dans le cas parabolique la solution de Lagrange peut \^etre \ep crite sous la forme suivante  $$
(q_1,q_2,q_3)=(q,\displaystyle\frac{q}{2},\displaystyle\frac{\sqrt{3}q}{2}),
\quad (p_1,p_2,p_3)=(p,\alpha
p+\displaystyle\frac{\beta}{q},\gamma p+\displaystyle\frac{\delta
}{q}), $$ o\`u $\alpha$, $\beta$, $\gamma$, $\delta$ sont les
constantes.

Pour $q(t)$, $p(t)$ on a  $$ q=P(w),\quad
p=\displaystyle\frac{w}{P(w)}, \quad P(w)=e_1w^2+e_2w+e_3,$$
 avec les constantes $e_1$, $e_2$, $e_3$.

 Nous avons donc  la param\ep trisation de la solution de Lagrange
$$ \Gamma = \{q_i(w),p_i(w) \mid w \in{\bf CP^1} \}. $$

A partir de cette solution particuli\el re, on peut calculer les {
\it \ep quations normales variationnelles}  (voir [5]) qui
prennent la forme d'un syst\el me fuchsien

$$  \displaystyle\frac{dx}{d\tau}=\left(
\displaystyle\frac{A}{\tau-\tau_0}+\displaystyle\frac{B}{\tau-\tau_1}+\displaystyle\frac{C}{\tau
-\tau_2}\right)x, \quad x\in {\bf C^4},  \eqno (2.1)
 $$
 o\`u $\tau_0$, $\tau_1$, $\tau_2 \in{\bf C}$ sont
les singularit\ep s et $A$, $B$, $C$ sont les $4\times 4$ matrices
constantes d\ep pendantes sur des masses $m_1$, $m_2$, $m_3$.

\vspace{0.5cm}

\noindent { \bf  3. Le groupe de monodromie des \ep quations
(2.1)}

 Soit $G$ le groupe de monodromie des \ep quations (2.1).

Nous d\ep signons  $T_i$, $T_{\infty}$, $i=0,1,2$  les g\ep n\ep
rateurs de $G$ correspondants respectivement aux groupes de
monodromie locaux autour les singularit\ep s $\tau_i$, $i=0,1,2$
et $\tau=\infty$. \vspace{0.5cm}

\noindent {\it LEMME 3.1} \\ \noindent a)  $T_0=I$ -- est l'\ep
l\ep ment neutre de $G$.

\noindent b) Les g\ep n\ep rateurs $T_1$, $T_2$ ont la m\^eme
forme de Jordan $$ \left(
\begin{array}{cccc}
1 & 1 & 0 & 0 \\ 0 & 1 & 0 & 0 \\ 0 & 0 & 1 & 1
\\ 0 & 0 & 0 & 1
\end{array}
\right).$$

 \noindent c) Le g\ep n\ep rateur $T_{\infty}$ a les valeurs propres suivantes   $$
{\mathrm Spectre}(T_{\infty})=\left\{ e^{2\pi i\lambda_1},\quad
e^{2\pi i\lambda_2},\quad e^{-2\pi i\lambda_1},\quad e^{-2\pi
i\lambda_2}\right\}, $$ o\`u $$
\lambda_1=\displaystyle\frac{3}{2}+\displaystyle\frac{1}{2}\sqrt{13+\sqrt{\theta}},\quad
\lambda_2=\displaystyle\frac{3}{2}+\displaystyle\frac{1}{2}\sqrt{13-\sqrt{\theta}},$$
et $$ \theta=144\left(1-\displaystyle
3\frac{S_2}{S_1^2}\right),\quad S_1=m_1+m_2+m_3,\quad
S_2=m_1m_2+m_3m_2+m_1m_3.$$

\vspace{0.5cm}

 \noindent {\it COROLLAIRE 3.2}   $$
T_1T_2=T_{\infty}^{-1}. $$

 \noindent {\it COROLLAIRE 3.3--} On aura toujours $$  {\mathrm Spectre}(T_{\infty} ) \neq \{1,1,1,1\}. $$

\noindent {\it Esquisse de d\ep monstration.--}
 Pour tout les singularit\ep s $\tau_i$, $\tau=\infty$, $i=0,1,2$ on
 peut calculer des solution formelles locales du syst\el me (2.1)
 et trouver donc les g\ep n\ep rateurs
de $G$ dans chaque point singulier.

\vspace{0.5cm}

 \noindent { \bf  4. Esquisse  de la d\ep monstration du th\ep or\el me 1.1 }

\vspace{0.5cm}

Notre d\ep monstration est inspir\ep e par Ziglin [7]. Supposons
que le probl\el me plan des trois corps (1.2) poss\el de deux
int\ep grales premi\el res  supplementaires m\ep romorphes  par
rapport aux variables $q_i$, $p_j$, $r_k$ et fonctionnellement
ind\ep pendantes au voisinage de la solution particuli\el re de
Lagrange $\Gamma$. Nous en d\ep duisons le lemme suivant

\vspace{0.5cm}
 \noindent{\it LEMME 4.1} (Ziglin [7])--Le groupe de monodromie $G$ des \ep quations normales
variationnelles (2.1) a deux invariants rationnels et
fonctionnellement ind\ep pendantes $J_1(x)$, $J_2(x)$.

\vspace{0.5cm}

 En vertu du lemme 3.1 on peut supposer sans perte
de g\ep n\ep ralit\' e que $$T_1= \left(
\begin{array}{cccc}
1 & 1 & 0 & 0 \\ 0 & 1 & 0 & 0 \\ 0 & 0 & 1 & 1
\\ 0 & 0 & 0 & 1
\end{array}
\right) =I+D,$$ o\`u  $$D=\left(
\begin{array}{cccc}
0& 1 & 0 & 0 \\ 0 & 0 & 0 & 0 \\ 0 & 0 & 0 & 1
\\ 0 & 0 & 0 & 0
\end{array}
\right). $$ On aura pour la matrice $T_2$ $$ T_2=I+R,$$ o\`u
$$R=V^{-1}D V=\left(
\begin{array}{cccc}
a_1 & a_2 & a_3 & a_4 \\ b_1 & b_2 & b_3 & b_4 \\ c_1 & c_2 & c_3
& c_4 \\ d_1 & d_2 & d_3 & d_4
\end{array}
\right),$$ avec une matrice $V$, $det(V)\neq 0$ et les constantes
inconnues $a_i$, $b_i$, $c_i$, $d_i\in {\bf C}.$

\vspace{0.5cm}

 \noindent{\it Remarque 4.2.} La matrice \'etant nilpotente,
 $${\mathrm Spectre}(R)=\{0,0,0,0\}. \eqno
 (4.1)$$

\vspace{0.5cm}

 \noindent{\it LEMME 4.3--}
Soit $J$ l'invariant rationnel de $G$, alors $$ \delta J=0, \quad
\Delta J =0,$$ o\`u $
\delta=x_2\displaystyle\frac{\partial}{\partial
x_1}+x_4\displaystyle\frac{\partial}{\partial x_3}, \quad \Delta=
\left(\sum_{i=1}^4a_ix_i\right)\displaystyle\frac{\partial}{\partial
x_1}+\left(\sum_{i=1}^4b_ix_i\right)\displaystyle\frac{\partial}{\partial
x_2}+\left(\sum_{i=1}^4c_ix_i\right)\displaystyle\frac{\partial}{\partial
x_3}+\left(\sum_{i=1}^4d_ix_i\right)\displaystyle\frac{\partial}{\partial
x_4}. $

\vspace{0.5cm}

 \noindent{\it D\ep monstration du lemme 4.3.--} Pour un arbitraire $n\in {\bf
 N }$  nous avons  $T_1^n=I+nD$. Par cons\ep quent
 $$
J\left(T_1^n x\right)=J(x+nDx)=J(x)+n\delta
J(x)+\sum\limits_{i=2}^{\infty} n^i r_i(x)=J(x),$$

Il r\ep sulte de l\`a  que $\delta J=0$. La preuve de l'identit\'
e $\Delta J=0$ s'effectue de la m\^eme mani\el re. \quad \quad
$\Box$

On obtient donc  $$ \delta J_i=\Delta J_i=0, \quad i=1,2.$$

Ces \ep quations nous permettent de trouver les restrictions sur
les param\el tres $a_i$, $b_i$, $c_i$, $d_i$ qui s'\ep crivent
sous la forme d'un syst\ep me d'\ep quations alg\ep briques $$
\begin{array}{lll}
P_1(a,b,c,d)=0,\\ \quad \quad \cdot \quad \cdot \quad \cdot
\\
 P_m(a,b,c,d)=0,
\end{array} \eqno (4.2)$$ avec les polyn\^omes $P_1\ldots P_m$, $m \in {\bf N}.$

Par ailleurs, la condition  4.1 nous donne un autre syst\el me
d'\ep quations sur les m\^emes \\ param\el tres   $$
\begin{array}{lll}
L_1(a,b,c,d)=0,\\ \quad \quad \cdot \quad \cdot \quad \cdot
\\
 L_n(a,b,c,d)=0,
\end{array} \eqno (4.3) $$ avec les polyn\^omes $L_1\ldots L_n$, $n \in {\bf N}.$

On peut v\ep rifier par un calcul direct, que  pour $a_i$, $b_i$,
$c_i$, $d_i$  satisfaisant aux syst\el mes (4.2) et (4.3) on aura
toujours $$ Spectr(T_1T_2)=\{1,1,1,1\},$$ d'o\`u, en utilisant le
corollaire 3.2 $$ Spectr(T_{\infty})=\{1,1,1,1\}.$$

Or cela est contraire au corollaire 3.3. Ceci ach\el ve la d\ep
monstration du th\ep or\el me 1.1.

\vspace{0.5cm}

 \noindent { \bf  Remerciements}

\vspace{0.5cm}

Je  remercie L. Gavrilov and V. Kozlov pour m'avoir propos\' e le
sujet de ce travail et ses suggestions. Je veux \ep galement
remercier J.--P. Ramis, J.--A. Weil et D. Bucher pour l'int\ep
r\^et qu'ils ont port\' e  \`a ce travail.

\begin{center}
 { \bf  R\ep f\ep rences bibliographiques}
\end{center}

\noindent [1] H. Bruns, {\it Ueberdie Integrale des vierk$\ddot
o$rper Problems}, Acta Math. 11, p. 25-96, 1887.

\noindent [2] J.L. Lagrange, {\it Oeuvres}. Vol. 6, 272-292,
Paris, 1873.

\noindent [3] J.J. Morales-Ruiz, J.P. Ramis, {\it Galosian
Obstructions to integrability of Hamiltonian Systems}, Preprint,
1998.

\noindent [4] H. Poincar\'e, {\it Les m\'ethodes novelles de la
m\'ecanique c\'eleste}, vol. 1, chap. 5, Paris,1892.

\noindent [5] A. Tsygvintsev, {\it On the variational equations of
three--body problem near Lagrangian solutions}, Pr\ep publication
150 du Laboratoire de Math\ep matiques E. Picard, Universit\' e
Toulouse III, 1999.

\noindent [6] E.T. Whittaker, {\it A Treatise on the Analytical
Dynamics of particles and Rigid Bodies}, Cambridge University
Press, New York, 1970.

\noindent [7] S.L. Ziglin, {\it Branching of solutions and
non-existence of first integrals in Hamiltonian Mechanics I},
Func. Anal. Appl. 16, 1982.

\end{document}